\documentclass[12pt]{article}

 \usepackage{amsmath}
\usepackage{amsfonts}
\usepackage{amssymb, amsfonts, amsmath, nicefrac}
\usepackage{xcolor}
\usepackage{subfig}

\newcommand{\cm}{{\sf m}}

\newcommand{\ci}{{\mathfrak i}}

\newtheorem{proposition}{Proposition}[section]

\newcommand{\be}{\begin{equation}}
\newcommand{\ee}{\end{equation}}

\def\e{\varepsilon}
\def\O{\Omega}

\def\bbr{{\Bbb{R}}} 

\def\bbc{{\Bbb{C}}}
\def\bbf{{\Bbb{F}}}

\begin{document}


\title{\bf Rank one tensor completion problem}

\author{\bf Mohit Singh, Alexander Shapiro and Rui Zhang}
\date{ }
\maketitle
\vspace{-1.1cm}
\begin{center}
         School of Industrial and Systems
Engineering,\\
                Georgia Institute of
Technology,\\
              Atlanta, Georgia 30332-0205.
\end{center}

\noindent
{\bf Abstract.}
 In this paper we consider the  rank-one tensor completion problem. We address the question of existence  and uniqueness of the  rank-one solution. In particular  we show that the global uniqueness over the field  of real numbers can be verified in a polynomial time. We give examples showing that there is an essential difference between  the question of global uniqueness  over the fields of  real and complex numbers. Finally we briefly discuss the rank-one  approximation problem for noisy observations.
\\



\setcounter{equation}{0}
\section{Introduction}
\label{sec:intr}

The  problem of  recovering  matrix of  low-rank  from a
few  observed entries is known as the
matrix   completion problem. There is a substantial body of literature on that problem,  e.g.  \cite{candes2010matrix,candrecht,DavenportRomberg16,Eldar2012,
hastie2015matrix,klopp2014noisy,Pimentel2016,sxz19}.  Tensor
completion is a natural   generalization  where the
goal is to recover a low-rank tensor from observations of   few of its  entries (e.g., \cite{Gandy, Kressner, Signoretto}).
The concept of tensor rank is considerably more involved then its matrix counterpart  (see, e.g.,  survey paper \cite{kolda2009tensor}). Although the low rank approximation problem has been
well studied for matrices, there is not much work on tensors.
From a generic point of view the question of  existence and uniqueness of low rank tensor decomposition   was investigated, e.g.,  in  \cite{chian2017} and references therein.

In this paper we consider the  rank-one tensor completion problem. In particular we address the question of uniqueness of   rank-one solution.
Necessary and sufficient conditions for uniqueness of rank-one {\em matrix } completion solution  are known and have a combinatorial flavor (cf., \cite{Kiral, sxz19}). On the other hand,
 tensor completion problem is considerably more involved and even rank-one  tensor completion   turns out to be nontrivial (cf., \cite{Kahle}).  Among other things we demonstrate  that  there is a significant difference between considering the problem over real and complex numbers.

\setcounter{equation}{0}
\section{Tensor completion}
\label{sec-tens}

Let us consider the following tensor completion problem over the fields of real or complex numbers.
Consider a three way tensor $X\in \bbf^{n_1\times n_2\times n_3}$ over the field of complex, $\bbf=\bbc$, or real $\bbf=\bbr$, numbers. Let $\O$ be a (nonempty)  set of indexes $(i,j,k)$, that is  $\O\subset [n_1]\times [n_2]  \times  [n_3]$,  and let  $Q_{ijk}$, $(i,j,k)\in \O$,  be the observed values (by $[n]$ we denote the set $\{1,...,n\}$). In the deterministic setting we want to find tensor $X=a\circ b\circ c$ of rank one such that $X_{ijk}=Q_{ijk}$, $(i,j,k)\in \O$. Recall that  $[a\circ b\circ c]_{ijk}=a_ib_jc_k$ with $a\in \bbf^{n_1}$, $b\in \bbf^{n_2}$, $c\in \bbf^{n_3}$. That is we want to find vectors $a,b,c$ such that
\begin{equation}\label{rank}
Q_{ijk}=a_ib_jc_k,\;(i,j,k)\in \O.
\end{equation}
We assume  throughout the paper that all observed values are nonzero, i.e.,
\begin{equation}
 Q_{ijk}\ne 0,\;(i,j,k)\in \O.
 \end{equation}
 This implies that the   components of vectors $a,b,c$  are nonzero.

The main goal of this paper is investigation of uniqueness of  the rank one solution. Note that vectors $a,b,c$  of   tensor  $a\circ b\circ c$ are defined up to change of scale $\lambda_1 a\circ \lambda_2 b\circ \lambda_3 c$,  where scalars
$\lambda_1   \lambda_2  \lambda_3 =1$.  Therefore we assume that $a_1=1$ and $b_1=1$.
 It is said that the solution  \eqref{rank} is {\em locally} unique if the representation \eqref{rank} is unique up to sufficiently small perturbations of vectors $a,b,c$.

By setting   $q_{ijk}:=\log |Q_{ijk}|$ we have the following system of linear equations associated with the representation \eqref{rank},
\begin{equation}\label{lin-1}
x_i+y_j+z_k=q_{ijk},\;(i,j,k)\in \O,
\end{equation}
with respect to variables $x_i:=\log |a_i|$, $y_j:=\log |b_j|$, $z_k:=\log |c_k|$.
Since  $a_1=1$ and $b_1=1$, and hence $x_1=0$ and $y_1=0$, we have that the system \eqref{lin-1} has $\cm=|\O|$ equations and $ n_1+ n_2+ n_3-2$ unknowns. Consider the following condition.
\begin{itemize}
  \item [(A)]
 Setting
$x_1=0$ and $y_1=0$,
 the homogeneous  linear system
 \begin{equation}\label{hom}
x_i+y_j+z_k=0,\;(i,j,k)\in \O,
\end{equation}
   with
$\cm=|\O|$ equations and $ n_1+ n_2+ n_3-2$ unknowns, is nondegenerate, i.e,
has only zero  solution.
\end{itemize}
Note that the  rank one solution \eqref{rank}  is   locally   unique iff the   rank one solution
\begin{equation}\label{absol}
 |Q_{ijk}|=|a_i| |b_j| |c_k|, \;(i,j,k)\in \O,
\end{equation}
 for the respective  absolute values is locally   unique. Therefore we have the following result.

\begin{proposition}
\label{pr-1}
The rank one solution \eqref{rank} is locally  unique if and only if condition {\rm (A)} holds.
\end{proposition}

The above characterization of local uniqueness is the same for the fields of real and complex numbers. Of course the local uniqueness is necessary for the corresponding global  uniqueness of the solution.
\begin{itemize}
  \item
Unless stated otherwise, we assume from now on  that condition  (A)  holds.
\end{itemize}
As we are going to show,  condition (A) is not sufficient for  the global  uniqueness.

 Let us consider the question of (global) uniqueness over complex numbers.
 We can write complex numbers $a_i,b_j,c_k$ in the following form  $a_i=|a_i|e^{\ci \alpha_i}$, $b_j=|b_j|e^{\ci  \beta_j}$,  $ |c_k|e^{\ci \gamma_k}$, for some $\alpha_i,\beta_j,\gamma_k\in [0,2\pi)$ where $\ci^2=-1$. Consider the following systems of linear equations, in unknowns $\alpha_i,\beta_j,\gamma_k$, $(i,j,k)\in [n_1]\times [n_2]\times [n_3]$,
\begin{equation}\label{comp-1}
\alpha_1=0,\;\beta_1=0,\; \alpha_i+\beta_j+\gamma_k=\sigma_{ijk},\;(i,j,k)\in \O,
\end{equation}
where $\sigma_{ijk}$  are either  0, $2\pi$ or $4\pi$, $(i,j,k)\in \O$.
Let $\bar{\alpha},\bar{\beta},\bar{\gamma}$ be a solution of \eqref{comp-1} for some $\sigma_{ijk}$. Then
\[
a_i b_j c_k=|a_i| |b_j| |c_k|e^{\ci(\alpha_i+\beta_j+\gamma_k)}=|a_i| |b_j| |c_k|e^{\ci(\alpha_i+\bar{\alpha}_i+\beta_j+\bar{\beta}_j+
\gamma_k+\bar{\gamma}_k)}.
\]
giving another solution $a_ie^{\ci\bar{\alpha}_i}$, $b_je^{\ci  \bar{\beta}_j}$,  $c_ke^{\ci \bar{\gamma}_k}$.
That is we have the following.
\begin{proposition}
\label{pr-unique}
 The rank one solution is {\em  not unique}, over the field of complex numbers,  if and only if
there exist $\alpha_i,\beta_j,\gamma_k\in (0,2\pi)$ which solve at least one   system of   equations    \eqref{comp-1} for $\sigma_{ijk}$  being either  $0$, $2\pi$ or $4\pi$, $(i,j,k)\in \O$.
\end{proposition}

Of course there are $3^{\cm}$ such systems of equations.
Note that the homogeneous counterpart of \eqref{comp-1} is the same as \eqref{hom}. Therefore (under condition (A))  if the system \eqref{comp-1} has a solution for some $\sigma_{ijk}$, it is unique. Note also that $a_i, b_j,c_k$ are determined by the respective $\alpha_i,\beta_j,\gamma_k$ up to the period of $2\pi$.

\setcounter{equation}{0}
\section{Rank one solution over real numbers}
\label{sec-real}

In this section, we give the criteria for existence and uniqueness of the rank one solution over real numbers. Consider the linear system of equations \eqref{lin-1} in unknowns $x_i,y_j,z_k$, $(i,j,k)\in [n_1]\times [n_2]\times [n_3]$, with $x_1=0$ and $y_1=0$.  If this system does not have a solution, then the corresponding rank one solution does not exist. Therefore assume that the system \eqref{lin-1} has a solution. Then  existence and uniqueness of the rank one solution comes to verification of the signs in the right hand side of \eqref{absol}.

Consider the finite field $GF(2)$ consisting of two elements 0 and 1 (this can be viewed  as integers with arithmetic operations modulo 2, where  0 representing the class of even integers and 1 represents the class of odd integers).
Consider the following system of equations over the field $GF(2)$:
\begin{equation}\label{modul}
 \e_i+\nu_j+\eta_k=c_{ijk}, \;(i,j,k)\in \O, \;\e_1=0,\;\nu_1=0,
\end{equation}
with unknowns $ \e_i,\nu_j,\eta_k\in GF(2)$,
where we set  $c_{ijk}=1$ if $Q_{ijk}>0$, and $c_{ijk}=0$ if $Q_{ijk}<0$. Then the rank one solution does exist iff the system \eqref{modul} has a solution. In the later case   the rank one  solution is unique iff \eqref{modul} has a unique solution. By running the  Gaussian  elimination procedure to solve \eqref{modul}  over the field $GF(2)$,  it is possible to check the above existence and uniqueness conditions. The running time of the algorithm is polynomial in the input size.

\setcounter{equation}{0}
\section{Numerical experiments and counterexamples}
\label{sec:numer}

In this section, by numerical examples, we show that for some $\Omega$ the global uniqueness over the field of real numbers does not imply the global uniqueness over the field of complex numbers (recall that it is assumed that  condition (A) is satisfied).
Of course  the solutions over the field of real numbers are also the solutions over the field of complex numbers.

In Tables \ref{Tensor1}, \ref{Tensor2}, \ref{Tensor3} and \ref{Tensor4}  tensors with partial observations are considered. From the left to the right, there are the first, second, and third slice of the tensor, with  ``$*$'' denoting  the missing value.
In all examples  the tensors are of size $3\times 3\times 3$, $|\Omega| = 7$. In table \ref{Tensor1}, \ref{Tensor2} , and \ref{Tensor3}, all observed values  equal to 1, i.e. $Q_{ijk} = 1,  \forall (i,j,k)\in\Omega$. Moreover, it  can be checked  that condition (A) is satisfied in these  examples. Recall that our goal is to find vectors $a,b,c$, such that $Q_{ijk} = a_ib_jc_k$, for $(i,j,k)\in\Omega$. Table \ref{Tensor1} gives  an example with only one  solution  over both fields  $\bbr$ and $\bbc$. Table \ref{Tensor2} gives an example with two solutions over $\bbr$.
Table \ref{Tensor3} gives  an example with only one solution over $\bbr$ but multiple solutions over $\bbc$. Table \ref{Tensor4} gives an example where there  exist solutions over $\bbc$,  but there are  no solutions  over $\bbr$. This can be verified by  a straightforward checking of  all relevant combinations.  The respective  solutions   are   provided.

  \begin{table}[!h]
    \centering
    \parbox{0.4\linewidth}{
    \centering
    \caption{Tensor with only one  solution}
    \label{Tensor1}
    \begin{tabular}{ccc|ccc|ccc}
    1&1&1&1&*&*&1&*&*\\
    1&*&*&*&*&*&*&*&*\\
    1&*&*&*&*&*&*&*&*
    \end{tabular}
    }
    \parbox{0.2\linewidth}{
      \centering
      \caption*{Solution}
      \begin{tabular}{c|c|c}
      $a$&$b$&$c$\\
      1&1&1\\
      1&1&1\\
      1&1&1
      \end{tabular}
    }
  \end{table}

     \begin{table}[!h]
      \centering
      \parbox{0.4\linewidth}{
      \caption{Tensor with two solutions over $\bbr$}
      \label{Tensor2}
      \begin{tabular}{ccc|ccc|ccc}
      1&1&*&1&*&*&*&*&1\\
      1&*&*&*&*&*&*&*&*\\
      *&*&1&*&*&*&1&*&*
      \end{tabular}
      }
      \centering
      \parbox{0.2\linewidth}{
        \centering
        \caption*{Solution 1}
        \begin{tabular}{c|c|c}
        $a$&$b$&$c$\\
        1&1&1\\
        1&1&1\\
        1&1&1
        \end{tabular}
      }
      \parbox{0.2\linewidth}{
        \centering
        \caption*{Solution 2}
        \begin{tabular}{c|c|c}
        $a$&$b$&$c$\\
        1&1&1\\
        1&1&1\\
        -1&-1&-1
        \end{tabular}
      }
    \end{table}

  \begin{table}[!h]
    \centering
    \caption{Tensor with   globally unique over $\bbr$,  but multiple   over $\bbc$ solutions}
    \label{Tensor3}
    \begin{tabular}{ccc|ccc|ccc}
    1&*&*&*&*&1&*&1&*\\
    1&*&*&*&*&*&*&*&*\\
    *&1&*&1&*&*&*&*&1
    \end{tabular}\\
    \vspace{0.2in}
    \centering
    \parbox{0.2\linewidth}{
      \centering
      \caption*{Solution 1}
      \begin{tabular}{c|c|c}
      $a$&$b$&$c$\\
      1&1&1\\
      1&1&1\\
      1&1&1
      \end{tabular}
    }
    \parbox{0.2\linewidth}{
      \centering
      \caption*{Solution 2}
      \begin{tabular}{c|c|c}
      $a$&$b$&$c$\\
      1&1&1\\
      1&$e^{\frac{4i\pi}{3}}$&$e^{\frac{4i\pi}{3}}$\\
      $e^{\frac{2i\pi}{3}}$&$e^{\frac{2i\pi}{3}}$&$e^{\frac{2i\pi}{3}}$
      \end{tabular}
    }
    \parbox{0.2\linewidth}{
      \centering
      \caption*{Solution 3}
      \begin{tabular}{c|c|c}
      $a$&$b$&$c$\\
      1&1&1\\
      1&$e^{\frac{2i\pi}{3}}$&$e^{\frac{2i\pi}{3}}$\\
      $e^{\frac{4i\pi}{3}}$&$e^{\frac{4i\pi}{3}}$&$e^{\frac{4i\pi}{3}}$
      \end{tabular}
    }
  \end{table}

     \begin{table}[!h]
      \centering
      \parbox{0.4\linewidth}{
      \caption{Tensor with solutions only over $\bbc$}
      \label{Tensor4}
      \begin{tabular}{ccc|ccc|ccc}
      -1&*&1&*&*&*&*&*&*\\
      1&*&*&*&*&*&*&1&*\\
      *&1&*&1&*&*&*&*&1
      \end{tabular}
      }
      \centering
      \parbox{0.2\linewidth}{
        \centering
        \caption*{Solution 1}
        \begin{tabular}{c|c|c}
        $a$&$b$&$c$\\
        1&1&-1\\
        -1&$e^{\frac{i\pi}{2}}$&$e^{\frac{3i\pi}{2}}$\\
        $e^{\frac{i\pi}{2}}$&-1&$e^{\frac{i\pi}{2}}$
        \end{tabular}
      }
      \parbox{0.2\linewidth}{
        \centering
        \caption*{Solution 2}
        \begin{tabular}{c|c|c}
        $a$&$b$&$c$\\
        1&1&-1\\
        -1&$e^{\frac{3i\pi}{2}}$&$e^{\frac{i\pi}{2}}$\\
        $e^{\frac{3i\pi}{2}}$&-1&$e^{\frac{3i\pi}{2}}$
        \end{tabular}
      }
    \end{table}

\setcounter{equation}{0}
\section{Rank one approximation of noisy observations}
\label{sec:noise}

Condition (A) can hold only if the number of unknowns is not bigger than the number of equations, i.e., if
$
  n_1+ n_2+ n_3-2\le \cm.
$
If the strict inequality  holds, i.e.,
\begin{equation}\label{ineq}
  n_1+ n_2+ n_3-2<\cm,
\end{equation}
then the  set of vectors $[Q]_{(i,j,k)\in \O}$, in the respective space of dimension $\cm$, has Lebesgue measure zero. That is, if values $Q_{ijk}$, $(i,j,k)\in \O$, are viewed as random variables with a (joint) continuous distribution, and the inequality \eqref{ineq} holds, then with probability one the one rank problem   does not have an exact  solution.

Suppose that the inequality \eqref{ineq} holds. As it was pointed above, in that case  the exact  rank one solution does not exist
with probability one. Therefore it makes sense to talk about approximate solutions. Let us consider the following model
\begin{equation}\label{model-1}
 Q_{ijk}=Q^*_{ijk}+\e_{ijk},\;(i,j,k)\in \O,
\end{equation}
where $Q^*_{ijk}=a_ib_jc_k$, $(i,j,k)\in \O$, allow exact one rank solution and
$\e_{ijk}$ are  viewed as noise or  disturbances (not necessary random). Suppose that  $Q^*_{ijk}>0$, $(i,j,k)\in \O$, and
$\e_{ijk}$ are small enough so that $Q_{ijk}>0$, $(i,j,k)\in \O$.
Let $q_{ijk}:=\log Q_{ijk}$ and consider
 the least squares problem
\begin{equation}\label{least}
 \min_{x,y,z}\sum_{(i,j,k)\in \O}(q_{ijk}-x_i-y_j-z_k)^2.
\end{equation}
Suppose that the condition (A) holds for $Q^*_{ijk}$, and let $\hat{x},\hat{y},\hat{z}$ be solutions of \eqref{least} with  $\hat{x}_1=\hat{y}_1=0$ (condition (A) ensures that these solutions are unique).

Note that
\[
q_{ijk}=\log Q^*_{ijk}+\log(1+\e_{ijk}/Q^*_{ijk})=q^*_{ijk}+\e^*_{ijk}+o(\e^*_{ijk}),
\]
where $q^*_{ijk}=\log Q^*_{ijk}$ and $\e_{ijk}^*= \e_{ijk}/Q^*_{ijk}$.
Therefore if the relative disturbances $\e^*_{ijk}$ are small,
then $\hat{a}_i=\exp(\hat{x}_i)$, $\hat{b}_j=\exp(\hat{y}_j)$, $\hat{c}_k=\exp(\hat{z}_k)$, could give good approximate solution for fitting one rank tensor to the observed values $Q_{ijk}$.

\subsection{Example}
In this example, the true tensor is of size $3\times 3\times 3$ and all the elements of it equal to 1, and $\varepsilon_{ijk}\sim \text{Uniform}(-0.2,0.2)$.
Noisy observations of the tensor are shown in table \ref{Tensor1Nsy}.
Table \ref{Tensor1NsySol} is the approximate solution from the least square problem, eq.(\ref{least}).
\begin{table}[!h]
  \centering
  \caption{Tensor with noisy observations}
  \label{Tensor1Nsy}
  \begin{tabular}{ccc|ccc|ccc}
  1.1718&1.1438&0.8739  &0.8193&0.8585&*    &0.9003&1.1636&*\\
  *&*&*                 &*&*&0.9160         &*&0.8386&1.0515\\
  0.8469&*&1.1119       &1.0942&0.8058&*    &0.9664&*&*
  \end{tabular}
\end{table}

\begin{table}[!h]
  \centering
  \caption{Rank-one approximate solution}
  \label{Tensor1NsySol}
  \begin{tabular}{ccc|ccc|ccc}
    1.0052&    1.0135&    1.0501&    0.8976&    0.9050&    0.9377&    0.9925&    1.0007&    1.0368\\
    0.9448&    0.9526&    0.9869&    0.8436&    0.8506&    0.8813&    0.9328&    0.9405&    0.9745\\
    0.9934&    1.0016&    1.0378&    0.8871&    0.8944&    0.9267&    0.9809&    0.9889&    1.0247
  \end{tabular}
  \vspace{0.2in}
    \centering
    \caption*{Solution}
    \begin{tabular}{c|c|c}
    $\hat a$&$\hat b$&$\hat c$\\
    1& 1&1.0052\\
    0.9399& 1.0082& 0.8976\\
    0.9883& 1.0446& 0.9925
  \end{tabular}
\end{table}

 \setcounter{equation}{0}
\section{Conclusion remarks}
\label{sec:conc}

It  could be mentioned that the approach   discussed in this paper can be extended in a straightforward way to an  analysis of rank one solutions of tensors of higher order. This also  could be  applied to   matrix completion,  of course a  two way tensor can be viewed as a matrix. For the rank one matrix completion problem the counterpart of the linear system \eqref{lin-1}
becomes
\begin{equation}\label{lin-mat}
x_i+y_j=q_{ij},\;(i,j)\in \O.
\end{equation}
By setting $x_1=0$  we have that the respective rank one solution is locally unique iff the homogeneous counterpart of equations  \eqref{lin-mat} is nondegenerate, i.e., has only zero solution. By
applying  the Gauss elimination procedure it is not difficult to see that
this condition is also necessary and sufficient for the global uniqueness and
is equivalent to conditions for uniqueness of  the rank one  matrix  completion solution derived  in  \cite{Kiral} and \cite{sxz19}.  On the other hand, as it was demonstrated in Section~\ref{sec:numer}, for the rank one tensor completion problem the situation is different - the local uniqueness does not imply the respective global uniqueness and conditions for global uniqueness over real and complex numbers are different.

\bibliographystyle{plain}
\bibliography{references}

\end{document}